\documentclass{amsart}
\usepackage{amssymb}
\usepackage{amscd}
\usepackage[latin1]{inputenc}
\usepackage{amsmath}
\usepackage{amsfonts}
\usepackage{latexsym}
\usepackage{verbatim}
\usepackage{epsfig}
\usepackage{calc}
\usepackage{color}
\usepackage{bbm}
\usepackage{cancel}

\newtheorem{theo}{Theorem}[section]
\newtheorem{theorem}[theo]{Theorem}
\newtheorem{prop}[theo]{Proposition}
\newtheorem{lemma}[theo]{Lemma}
\newtheorem{corol}[theo]{Corollary}

\theoremstyle{definition}

\newtheorem{schritt}{Step}

\newtheorem{claim}{Claim}

\newtheorem{questions}[theo]{Questions}

\theoremstyle{remark}
\newtheorem{rem}[theo]{Remark}

\newtheorem*{ack}{Acknowledgment}

\newcounter{help}

\newsavebox{\loopl}

\newcommand{\ie}[0]{i.e.\ }
\newcommand{\eg}[0]{e.g.\ }

\newcommand{\resp}[0]{resp.\ }

\newcommand{\sft}{\textnormal{SFT}}

\newcommand{\ufp}{\textnormal{UFP}}

\newcommand{\abs}[1]{\left\lvert#1\right\rvert}
\newcommand{\norm}[1]{\left\lVert#1\right\rVert}

\newcommand{\zz}{{\mathbb Z}}
\newcommand{\nz}{{\mathbb N}}

\newcommand{\idop}{{\mathbbm{1}}}

\newcommand{\zt}{{\zz^2}}
\newcommand{\zd}{{\zz^d}}

\newcommand{\vj}{{\vec{\jmath}}}
\newcommand{\vi}{{\vec{\imath}}}
\newcommand{\vu}{{\vec{u}}}
\newcommand{\vv}{{\vec{v}}}
\newcommand{\vw}{{\vec{w}}}

\newcommand{\ag}{\ensuremath{{\mathcal A}} }

\newcommand{\cg}{\ensuremath{{\mathcal C}} }

\newcommand{\lang}{\ensuremath{{\mathcal L}} }
\newcommand{\mg}{\ensuremath{{\mathcal M}} }
\newcommand{\pg}{\ensuremath{{\mathcal P}} }

\newcommand{\seq}[1]{\left< #1 \right>}

\newcommand{\dcup}{\ensuremath{\,{\mathbin{\dot\cup}}\, }}

\newcommand{\wel}{W^{\textnormal{el}}}
\newcommand{\htop}{h_{\textnormal{top}}}

\makeatletter
\definecolor{gray}{rgb}{.8, .8, .8}
\definecolor{darkgray}{rgb}{.5, .5, .5}
\newcommand{\const}{\mathop{\operator@font const}\nolimits}
\newcommand{\card}{\mathop{\operator@font card}\nolimits}
\newcommand{\diam}{\mathop{\operator@font diam}\nolimits}
\newcommand{\id}{\mathop{\operator@font Id}\nolimits}
\newcommand{\orb}{\mathop{\operator@font Orb}\nolimits}
\newcommand{\fix}{\mathop{\operator@font Fix}\nolimits}
\newcommand{\erz}{\mathop{\operator@font span}\nolimits}
\newcommand{\supp}{\mathop{\operator@font supp}\nolimits}
\makeatother

\begin{document}
\title {Projectional entropy and the electrical wire shift -- Preliminary version}
\begin{abstract}
In this paper we present an extendible, block gluing $\zz^3$ shift of finite type $\wel$ in which the topological entropy equals the $L$-projectional entropy for a two-dimensional sublattice $L\subsetneq\zz^3$, even so $\wel$ is not a full $\zz$ extension of $\wel_L$. In particular this example shows that Theorem 4.1 of \cite{jkm} does not generalize to $r$-dimensional sublattices $L$ for $r>1$.

Nevertheless we are able to reprove and extend the result about one-dimensional sublattices for general $\zd$ shifts under the same mixing assumption as in \cite{jkm} and by posing a stronger mixing condition we also obtain the corresponding statement for higher-dimensional sublattices.
\end{abstract}

\date{}
\author{Michael Schraudner}
\address{Michael Schraudner\\
Centro de Modelamiento Matematico\\ 
Universidad de Chile\\
Av. Blanco Encalada 2120, Piso 7\\
Santiago de Chile}
\email{mschraudner@dim.uchile.cl}
\urladdr{www.math.uni-heidelberg.de/studinfo/schraudner}
\thanks{The author was supported by FONDECYT project 3080008.}
\renewcommand{\subjclassname}{MSC 2000}
\maketitle

\section{Preliminaries}\label{preliminariessection}

The purpose of this paper is to construct an example shift of finite type showing that Theorem 4.1 in \cite{jkm} is no longer valid -- contrary to an assertion made by the authors -- in the case of higher-dimensional sublattices. Nonetheless we get the claimed result even for general $\zd$ shifts instead of extendible shifts of finite type by imposing a stronger mixing condition. The statement about one-dimensional sublattices also holds for general $\zd$ shifts even in the original setting, \ie assuming only the weaker mixing property of \cite{jkm}.

We assume a basic familiarity with (multidimensional) symbolic dynamics, thus here we just fix some notation.

Every finite alphabet $\ag$ gives rise to a {\it $d$-dimensional full shift} $\ag^\zd$ ($d\in\nz$), a space equipped with the product topology of the discrete topology on $\ag$ which supports a natural $\zd$ (shift) action $\sigma: \zd\times\ag^\zd\rightarrow\ag^\zd$ given by translation $(\sigma_\vi(x))_\vj=(\sigma(\vi,x))_\vj:=x_{\vi+\vj}$ for all $\vi,\vj\in\zd$, $x\in\ag^\zd$.

Any closed $\sigma$-invariant subset $X\subseteq\ag^\zd$ together with the restricted shift action $\sigma|_{\zd\times X}$ constitutes a {\it $\zd$ shift}. If $X$ can be defined using a finite set $\pg\subseteq\ag^F$ of allowed patterns on some finite non-empty shape $F\subsetneq\zd$, so that $X=\{x\in\ag^\zd\mid \forall\,\vi\in\zd:\;x|_{\vi+F}\in\pg\}$, it is called a {\it $\zd$ shift of finite type} (\sft).

In the following we will use $\lang_F(X)$ to denote the set of patterns $\{x|_F\mid x\in X\}$ appearing in $X$ on some fixed finite subset $F\subsetneq\zd$ of coordinates. The {\it language} $\lang(X)$ of a $\zd$ shift $X$ consisting of all finite patterns that occur as subwords of elements of $X$ then is the union of $\lang_F(X)$ over all $F\subsetneq\zd$ finite.

The {\it topological entropy} of a $\zd$ shift $X$ measures the exponential growth rate of patterns and is defined in complete analogy to the one-dimensional setting as
\[
\htop(X):=\lim_{n\rightarrow\infty}\frac{\log\abs{\lang_{C_n}(X)}}{\abs{C_n}}
\]
where $C_n:=\{\vi\in\zd\mid \norm{\vi\,}_\infty\leq n\}$. Often we will write $h(X)$ instead of $\htop(X)$.

Following \cite{jkm} we now define sublattices, projectional entropy and degeneracy of $\zd$ shifts: For $d\in\nz$ and $1\leq r<d$ let ${\mathcal U}=\{\vu^{(1)},\ldots,\vu^{(r)}\},{\mathcal V}=\{\vv^{(1)},\ldots,\vv^{(d-r)}\}\subsetneq\zd$ be two disjoint sets of integer vectors such that ${\mathcal U}\dcup{\mathcal V}$ is a linearly independent set whose integer span $\erz_\zz({\mathcal U}\dcup{\mathcal V})=\seq{\vu^{(1)},\ldots,\vu^{(r)},\vv^{(1)},\ldots,\vv^{(d-r)}}_\zz$ equals $\zd$. Then $L:=\erz_\zz({\mathcal U})=\seq{\vu^{(1)},\ldots,\vu^{(r)}}_\zz\subsetneq\zd$ is called an {\it $r$-dimensional sublattice} of $\zd$. (Using $\vu^{(1)},\ldots,\vu^{(r)}$ as generators, $L$ is isomorphic to $\zz^r$.) The set $L':=\erz_\zz({\mathcal V})=\seq{\vv^{(1)},\ldots,\vv^{(d-r)}}_\zz$ constitutes a complementary $(d-r)$-dimensional sublattice.

Let $X$ be some $\zd$ shift and $L$ be any $r$-dimensional sublattice ($1\leq r<d$). By projecting points of $X$ onto $L$ we obtain a $\zz^r$ shift $X_L:=\{x|_L\mid x\in X\}$ on which the $\zz^r$ shift action is given as $\sigma|_{L\times X_L}$. Now the {\it $L$-projectional entropy} of $X$ is the topological entropy of the $\zz^r$ shift $X_L$ and we denote this quantity by $h_L(X):=\htop(X_L)$.

Given a $\zd$ shift $X$ and some $r$-dimensional sublattice $L$ ($1\leq r<d$) we form a new $\zd$ shift by taking the cartesian product of $X_L$ with itself along some complementary sublattice $L'$: $(X_L)^{\zz^{d-r}}:=\prod_{L'} X_L=\{(x^{(\vw)}\in X_L)_{\vw\in L'}\}$. For every $\vu\in L$, $\vv\in L'$ the symbol at coordinate $\vu+\vv\in\zd$ in the point $(x^{(\vw)}\in X_L)_{\vw\in L'}$ is given by $(x^{(\vv)})_\vu$. Note that by shift-invariance this construction is independent of the complementary sublattice $L'$ we chose. Obviously $X\subseteq (X_L)^{\zz^{d-r}}$. In the case that $X=(X_L)^{\zz^{d-r}}$, \ie $X$ is a full $\zz^{d-r}$ extension of $X_L$, we call $X$ {\it degenerate (with respect to $L$)}.

Fixing $\vu,\vv\in\zd$ the finite set $B:=\{\vi\in\zd\mid \forall\, 1\leq k\leq d:\ \vu_k\leq\vi_k\leq\vv_k\}$ is called a {\it (rectangular/cuboid) block} and we will use the notation $B=[\vu,\vv\,]$ to denote this set. Moreover we set $\vec{\idop}\in\zd$ to be the vector with all its components equal to $1$, thus for $n\in\nz_0$ we have $[-n\vec{\idop},n\vec{\idop}]=\{\vi\in\zd\mid \norm{\vi\,}_\infty\leq n\}$.

We say a $\zd$ \sft\ $X$ is {\it extendible}, if for any block $B=[\vv,\vw]\subsetneq\zd$ every allowed configuration $P\in\ag^B$ actually is in $\lang_B(X)$, \ie every locally valid pattern $P$ on a block $B$ can be extended to a point in $X$.

We finish this section by recalling three uniform mixing properties. The first one considering only pairs of cubes $[-n\vec{\idop},n\vec{\idop}]$ comes from \cite{jkm}, whereas the second, more homogeneous looking one that takes into account arbitrary blocks was introduced in \cite{bps}.

A $\zd$ shift $X$ is called {\it block strongly irreducible} \cite{jkm} if there exists a constant $s\in\nz$ such that whenever $\vi,\vj\in\zd$ and $m,n\in\nz$ satisfy that the distance between the blocks $\vi+[-m\vec{\idop},m\vec{\idop}]$ and $\vj+[-n\vec{\idop},n\vec{\idop}]$ (with respect to the maximum-metric on $\zd$) is larger than $s$ any two patterns $P_1\in\lang_{[-m\vec{\idop},m\vec{\idop}]}(X),\; P_2\in\lang_{[-n\vec{\idop},n\vec{\idop}]}(X)$ can be put together, \ie there exists a point $x\in X$ with $x|_{\vi+[-m\vec{\idop},m\vec{\idop}]}=P_1$ and $x|_{\vj+[-n\vec{\idop},n\vec{\idop}]}=P_2$.

A $\zd$ shift $X$ is called {\it block gluing} \cite{bps} if there exists a constant $g\in\nz_0$ such that whenever the two blocks $B_1=[\vu^{(1)},\vv^{(1)}],\; B_2=[\vu^{(2)},\vv^{(2)}]\subsetneq\zd$ have a distance larger than $g$ any pair of patterns $P_1\in\lang_{B_1}(X),\;P_2\in\lang_{B_2}(X)$ can be put together, \ie there exists a point $x\in X$ with $x|_{B_1}=P_1$ and $x|_{B_2}=P_2$.

In fact it is not hard to show that even though defined seemingly different, both notions actually coincide and we will use them interchangeably.

\begin{lemma}
A $\zd$ shift is block strongly irreducible if and only if it is block gluing.
\end{lemma}

The third mixing property -- putting a strictly stronger condition on $X$ -- that will be used in our results is the {\it uniform filling property} (\ufp) introduced in \cite{rs}.

A $\zd$ shift $X$ has the \ufp\ \cite{rs} if there exists a filling length $l\in\nz_0$ such that whenever we take a point $y\in X$ and a pattern $P\in\lang_B(X)$ on some block $B=[\vu,\vv]\subsetneq\zd$ there exists a point $x\in X$ with $x|_B=P$ and $x|_{\zd\setminus[\vu-l\vec{\idop},\vv+l\vec{\idop}]}=y|_{\zd\setminus[\vu-l\vec{\idop},\vv+l\vec{\idop}]}$.

\section{Main results on projectional entropy}\label{mainresultssection}

In \cite{jkm} the following theorem about $2$-dimensional \sft s was proved and it was claimed \cite[page 250]{jkm} that the same result holds for general dimensions $d>2$ and any sublattice $L\subsetneq\zd$.

\begin{theo}[{\cite[Theorem 4.1]{jkm}}]\label{jkmtheorem}
Let $(X,\sigma)$ be an extendible, block strongly irreducible $\zt$ \sft\ and $L\subsetneq\zt$ a $1$-dimensional sublattice. Then $\htop(X)=h_L(X)$ if and only if $X=(X_L)^{\zz}$.
\end{theo}

However, to ensure a certain property of the projectional shift $X_L$, the primordial proof uses Corollary 4.4.9 from \cite{lm}, which is valid for one-dimensional shifts only. Hence the generalization to higher-dimensional sublattices is not obvious and in fact is not true in general: In Section \ref{Welconstructionsection} we construct a $\zz^3$ \sft\ named the ``electrical wire shift'' that provides an example for which the above theorem does not hold in the case of a two-dimensional sublattice.

\begin{prop}\label{WelCounterexampleprop}
The electrical wire shift $\wel$ (see Section \ref{Welconstructionsection}) is an extendible, block strongly irreducible $\zz^3$ \sft\ such that $\htop(\wel)=h_L(\wel)$ for the sublattice $L:=\seq{\vec{e}_1,\vec{e}_2}_\zz\subsetneq\zz^3$. Nevertheless $\wel\subsetneq(\wel_L)^\zz$.
\end{prop}

We defer the proof of Proposition \ref{WelCounterexampleprop} to Section \ref{Welconstructionsection} where all necessary properties of $\wel$ are shown during its construction.

Existence of the electrical wire shift shows that the argument of Johnson, Kass and Madden used in the proof of their Theorem 4.1 does NOT extend to sublattices of dimension $2$ or higher without additional conditions. One possible such assumption is a stronger mixing property like the \ufp. Imposing this we obtain the following general theorem. (Note that we do not assume our $\zd$ shift to be extendible or of finite type.)

\begin{theorem}\label{UFPTheorem}
Let $X$ be a $\zd$ shift with the uniform filling property and let $L\subsetneq\zd$ be any $r$-dimensional sublattice ($1\leq r<d$), then $\htop(X)=h_L(X)$ if and only if $X=(X_L)^{\zz^{d-r}}$.
\end{theorem}

\begin{proof}
The implication ``$X=(X_L)^{\zz^{d-r}}$, thus $\htop(X)=h_L(X)$'' is trivial.

For the converse we use that if $X$ has the \ufp, the same is true for the cartesian product $(X_L)^{\zz^{d-r}}$. (Observe that we do not claim $X_L$ seen as a $\zz^r$ shift has to have the \ufp\ though.)

To show this, suppose $X$ has filling length $l\in\nz$ and $L'\subsetneq\zd$ is some $(d-r)$-dimensional sublattice complementary to $L$. For any point $y=\bigl(y^{(\vw)}\in X_L\bigr)_{\vw\in L'}\in (X_L)^{\zz^{d-r}}$ we can take a family $\bigl(\widetilde{y}^{(\vw)}\in X\bigr)_{\vw\in L'}$ of preimages under the projection onto $L$ such that for all $\vw\in L'$ we have $\widetilde{y}^{(\vw)}|_{\vw+L}=y^{(\vw)}=y|_{\vw+L}$. Similarily for any finite pattern $P\in\lang_B\bigl((X_L)^{\zz^{d-r}}\bigr)$ on some block $B=[\vu,\vv]$ ($\vu,\vv\in\zd$) we have a point $z=\bigl(z^{(\vw)}\in X_L\bigr)_{\vw\in L'}\in(X_L)^{\zz^{d-r}}$ realizing $P$, \ie $z|_B=P$. Again we have a family of preimages $\bigl(\widetilde{z}^{(\vw)}\in X\bigr)_{\vw\in L'}$ with $\widetilde{z}^{(\vw)}|_{(\vw+L)\cap B}=z|_{(\vw+L)\cap B}$. Applying the uniform filling property in $X$ for every pair $\widetilde{y}^{(\vw)},\;\widetilde{z}^{(\vw)}|_B$ with $\vw\in L'$ we get a family of points $\bigl(\widetilde{x}^{(\vw)}\in X\bigr)_{\vw\in L'}$ such that $\widetilde{x}^{(\vw)}|_B=\widetilde{z}^{(\vw)}|_B$ and $\widetilde{x}^{(\vw)}|_{\zd\setminus[\vu-l\vec{\idop},\vv+l\vec{\idop}]}=\widetilde{y}^{(\vw)}|_{\zd\setminus[\vu-l\vec{\idop},\vv+l\vec{\idop}]}$. Now we can project each $\widetilde{x}^{(\vw)}$ onto the coordinates in $\vw+L$ to obtain a point $x^{(\vw)}\in X_L$. (More precisely we first shift $\widetilde{x}^{(\vw)}$ by $\vw\in L'\subsetneq\zd$, then project $\sigma_{\vw}(\widetilde{x}^{(\vw)})\in X$ onto $L$ and define $x^{(\vw)}$ as the image $\sigma_{\vw}(\widetilde{x}^{(\vw)})|_L\in X_L$.) Putting all the $x^{(\vw)}$ together we finally form a point $x=\bigl(x^{(\vw)}\in X_L\bigr)_{\vw\in L'}\in (X_L)^{\zz^{d-r}}$ and we can easily check that $x|_B=P$ and $x|_{\zd\setminus[\vu-l\vec{\idop},\vv+l\vec{\idop}]}=y|_{\zd\setminus[\vu-l\vec{\idop},\vv+l\vec{\idop}]}$. Thus $l$ is also a filling length for $(X_L)^{\zz^{d-r}}$.

It is known that the uniform filling property implies entropy minimality. For a detailed proof of this technical result see Lemma \ref{UFPentropyminlemma}. Hence $(X_L)^{\zz^{d-r}}$ having the \ufp\ forces a strict entropy inequality $\htop(X)<\htop\bigl((X_L)^{\zz^{d-r}}\bigr)=\htop(X_L)=h_L(X)$ whenever $X\subsetneq(X_L)^{\zz^{d-r}}$ is a proper subsystem. This finishes our proof.
\end{proof}

\begin{rem}
In \cite[Appendix C]{bps} we defined the meandering streams shift and we proved it to be a non entropy minimal $\zt$ \sft\ which is corner gluing in each of the $4$ corners (in $NE$-, $NW$-, $SE$-, $SW$-direction). Using a similar construction as for the electrical wire shift, this example shows that we can not weaken the mixing assumption in Theorem \ref{UFPTheorem} to 4-corner-gluing.
\end{rem}

In the case of $1$-dimensional sublattices we get the corresponding generalization of Theorem \ref{jkmtheorem} for $\zd$ shifts even assuming only the weaker (original) mixing assumption used in \cite{jkm}. The proof of the following theorem was found in collaboration with R.\ Pavlov \cite{ppc}.

\begin{theorem}\label{blockTheorem}
Let $X$ be a block gluing $\zd$ shift and let $L\subsetneq\zd$ be any $1$-dimensional sublattice, then $\htop(X)=h_L(X)$ if and only if $X=(X_L)^{\zz^{d-1}}$.
\end{theorem}

\begin{proof}
Again one of the implications is trivial. For the converse let $X$ be block gluing with gluing constant $g\in\nz_0$ and let $L=\seq{\vw}_\zz$ be generated by $\vw\in\zd$. We claim that $X_L$ is block gluing as well.

For this let $P_1\in\lang_{B_1}(X_L)$ and $P_2\in\lang_{B_2}(X_L)$ be any two finite words on $L$-intervals $B_1=[u^{(1)},v^{(1)}]:=\{j\vw\mid u^{(1)}\leq j\leq v^{(1)}\},B_2=[u^{(2)},v^{(2)}]\subsetneq L$. Thus we have points $y,z\in X_L$ with $y|_{B_1}=P_1$, $z|_{B_2}=P_2$ and taking preimages of those gives $\widetilde{y},\widetilde{z}\in X$ such that $\widetilde{y}|_L=y$, $\widetilde{z}|_L=z$. Define $\zd$-blocks $\widetilde{B}_i:=[u^{(i)}\vw,v^{(i)}\vw]\subsetneq\zd$ ($i=1,2$). If we assume $v^{(1)}+g<u^{(2)}$ those blocks $\widetilde{B}_1,\widetilde{B}_2$ are at least a distance $g+1$ apart from each other. Since $X$ is block gluing there exists a point $\widetilde{x}\in X$ with $\widetilde{x}|_{\widetilde{B}_1}=\widetilde{y}|_{\widetilde{B}_1}$ and $\widetilde{x}|_{\widetilde{B}_2}=\widetilde{z}|_{\widetilde{B}_2}$. Projecting $\widetilde{x}$ onto $L$ we get a point $x:=\widetilde{x}|_L\in X_L$ which realizes the two patterns $P_1,P_2$ exactly at $B_1,B_2$ and hence $g$ is also a gluing constant for $X_L$.

Now recall that for $\zz$ shifts to be block gluing is equivalent to having the specification property -- for a definition see \cite[Section 21]{dgs}. Next we exploit that expansive dynamical systems with specification are intrinsically ergodic \cite[Theorem 22.15]{dgs}. Hence $X_L$ carries a unique measure of maximal entropy usually called the Bowen measure $\mu_B\in\mg_{\textnormal{Max}}(X_L)$ which has full support \cite[Theorem 22.10, Proposition 22.17]{dgs}. Then the product measure $\widetilde{\mu}:=\mu_B^{\zz^{d-1}}\in \mg_{\textnormal{Max}}\bigl((X_L)^{\zz^{d-1}}\bigr)$ is the unique maximal measure (ref ?) on the cartesian product $(X_L)^{\zz^{d-1}}$. Since we assume $\htop(X)=h_L(X)=\htop\bigl((X_L)^{\zz^{d-1}}\bigr)$, every measure $\nu\in\mg_{\textnormal{Max}}(X)$ of maximal entropy for $X$ is at the same time a measure of maximal entropy for $(X_L)^{\zz^{d-1}}$. Therefore $\nu=\widetilde{\mu}$ and $(X_L)^{\zz^{d-1}}=\supp(\widetilde{\mu})=\supp(\nu)\subseteq X$ implies $X=(X_L)^{\zz^{d-1}}$ as claimed.
\end{proof}

\begin{rem}
Observe that for sublattices $L\subsetneq\zd$ of dimension $r\geq2$ the $\zz^r$ shift $X_L$ in general does not inherit the block gluing property from $X$. For one such example use the electrical wire shift $\wel$ and the $2$-dimensional sublattice $L$ generated by $\vec{e}_1$ and $\vec{e}_1+\vec{e}_2$. Choosing $N\in\nz$, a pattern of blanks along the diagonal $B_1:=\{n(\vec{e}_1+\vec{e}_2)\mid 0\leq n\leq N\}\subsetneq L$ (which is a rectangular block in the $\zt$ shift $\wel_L$) forces blanks on all of $\{m\vec{e}_1+n\vec{e}_2\mid 0\leq m,n\leq N\}\subsetneq L$. Therefore it never can be put together with a non blank at $B_2=\{N\vec{e}_1\}$, even though the distance between the blocks $B_1$ and $B_2$ equals $N$ and thus can be made larger than any given gluing constant. So in the proof of Theorem \ref{blockTheorem} we implicitly had to use the special geometry of a one-dimensional lattice and this is the reason why the result is necessarily different from the one obtained for a higher dimensional sublattice.
\end{rem}

The following technical fact about the topological entropy of subsystems of $\zd$ shifts having the \ufp\ seems to be known at least for \sft s \cite{qs}, though we are not aware of an explicit demonstration in the literature. For completeness we include a proof which also gives the result in the case of general $\zd$ shifts.

\begin{lemma}\label{UFPentropyminlemma}
Every $\zd$ shift $X$ having the uniform filling property is entropy minimal, \ie every non-empty proper subsystem of $X$ has strictly smaller (topological) entropy.
\end{lemma}

\begin{proof}
Let $Y\subsetneq X$ be a proper subsystem of $X$. Hence there exists a pattern $P\in\lang(X)\setminus\lang(Y)$, say of shape $B=[\vec{\idop},n\vec{\idop}]$ for some $n\in\nz$. Assume $X$ has the \ufp\ with a filling length $l\in\nz$ and put $\widetilde{B}:=[(1-l)\vec{\idop},(n+l)\vec{\idop}]$.

We prove the following bound on the number of valid $Y$-patterns in comparison to the number of valid $X$-patterns on a large hypercube $C(N):=[\vec{\idop},N(n+l)\vec{\idop}]$ ($N\in\nz$).
\begin{equation}\label{patternbound}\tag{PB}
\abs{\lang_{C(N)}(Y)}\leq \bigl(1-\abs{\lang_{\widetilde{B}}(X)}^{-1}\bigr)^{N^d}\cdot \abs{\lang_{C(N)}(X)}
\end{equation}

For this let $J=\{\vj_i\mid 1\leq i\leq N^d\}:=\{\vj\in \vec{\idop}+(l+n){\nz_0}^d\mid\norm{\vj\,}_\infty\leq N(l+n)\}$ and for every subset $I\subseteq J$ define
\begin{align*}
\lang_{C(N)}^I(X):=
\{x|_{C(N)}\mid x\in X\wedge\forall\,\vj\in J\setminus I:\,x_{\vj+B}\neq P\}\ .
\end{align*}
Note that $\lang_{C(N)}^\emptyset(X)=\{x|_{C(N)}\mid x\in X\wedge\forall\,\vj\in J:\,x_{\vj+B}\neq P\}$ and $\lang_{C(N)}^J(X)=\lang_{C(N)}(X)$. Now for every $\vi\in I\subseteq J$ we have
\begin{align*}
\bigl|{\lang_{C(N)}^{I}(X)}\bigr| &\stackrel{{\ufp}}{\leq}\bigl|{\lang_{\vi+\widetilde{B}}(X)}\bigr|\cdot \bigl|{\bigl\{x|_{C(N)}\bigm\vert x\in X\wedge x|_{C(N)}\in\lang_{C(N)}^{I}(X)\wedge x|_{\vi+B}=P\bigr\}}\bigr|
\end{align*}
and using $\lang_{\vi+\widetilde{B}}(X)=\lang_{\widetilde{B}}(X)$ we get a lower bound
\begin{multline}\label{lowerbound}\tag{*}
\bigl|{\bigl\{x|_{C(N)}\bigm\vert x\in X\wedge x|_{C(N)}\in\lang_{C(N)}^{I}(X)\wedge x|_{\vi+B}=P\bigr\}}\bigr|\\
\geq\abs{\lang_{\widetilde{B}}(X)}^{-1}\cdot\bigl|{\lang_{C(N)}^{I}(X)}\bigr|
\end{multline}

Now we can estimate the number of valid $Y$-patterns on $C(N)$.
\begin{align*}
\abs{\lang_{C(N)}(Y)}&\leq\bigl|{\lang_{C(N)}^\emptyset(X)}\bigr|= \bigl|{\lang_{C(N)}^{\{\vj_1\}}(X)\setminus\bigl\{x|_{C(N)}\bigm\vert x\in X\wedge x|_{\vj_1+B}=P\bigr\}}\bigr|\\
&=\bigl|{\lang_{C(N)}^{\{\vj_1\}}(X)}\bigr|-\bigl|{\bigl\{x|_{C(N)}\bigm\vert x\in X\wedge x|_{C(N)}\in \lang_{C(N)}^{\{\vj_1\}}(X)\wedge x|_{\vj_1+B}=P\bigr\}}\bigr|\\
&\stackrel{\eqref{lowerbound}}{\leq}
\bigl(1-\abs{\lang_{\widetilde{B}}(X)}^{-1}\bigr)\cdot\bigl|{\lang_{C(N)}^{\{\vj_1\}}(X)}\bigr|\leq \bigl(1-\abs{\lang_{\widetilde{B}}(X)}^{-1}\bigr)^2\cdot\bigl|{\lang_{C(N)}^{\{\vj_1,\vj_2\}}(X)}\bigr|\\
&\leq\ldots\leq\bigl(1-\abs{\lang_{\widetilde{B}}(X)}^{-1}\bigr)^{\abs{J}}\cdot \bigl|{\lang_{C(N)}^{J}(X)}\bigr|\\
&=\bigl(1-\abs{\lang_{\widetilde{B}}(X)}^{-1}\bigr)^{N^d}\cdot\abs{\lang_{C(N)}(X)}
\end{align*}

This proves \eqref{patternbound} and putting this bound into the definition of topological entropy yields
\begin{align*}
\htop(Y)
&\leq\lim_{N\rightarrow\infty}\frac{{N^d}\log{\bigl(1-\abs{\lang_{\widetilde{B}}(X)}^{-1}\bigr)}+ \log{\abs{\lang_{C(N)}(X)}}}{N^d(l+n)^d}\\
&=\frac{1}{(l+n)^d}\underbrace{\log{\bigl(1-\abs{\lang_{\widetilde{B}}(X)}^{-1}\bigr)}}_{<0}+\htop(X)<\htop(X).
\end{align*}
Since $Y\subsetneq X$ was arbitrary, $X$ indeed is entropy minimal.
\end{proof}

\section{The construction of the electrical wire shift $\wel$}\label{Welconstructionsection}

We build $\wel$ in three steps: First we construct an extendible, block gluing $\zt$ \sft\ $W$ which models a system of straight wires running in a $\zt$ plane. A wire may branch into multiple subwires and those can unify again. However there are neither electrical sources nor consumers (sinks), thus once present, a wire has to go on forever without starting or ending at a certain coordinate in $\zt$. In a second step we slightly modify $W$ keeping the extendability and the block gluing property by independently replacing the occurences of a particular symbol with elements of a set of $k\geq 2$ distinct but completely interchangeable copies of this symbol. Doing this we get a family of new $\zt$ \sft s $\widetilde{W}_k$ which have larger entropy and in particular are no longer entropy minimal. Finally we put together $\zt$ configurations of $\widetilde{W}_2$ building our electrical wire shift $\wel\subsetneq\widetilde{W}_2^\zz$. Again $\wel$ will be extendible and block gluing. Points of $\wel$ projected onto the two-dimensional sublattice $L:=\seq{\vec{e}_1,\vec{e}_2}_\zz\subsetneq\zz^3$ will look like arbitrary $\zt$ configurations of $\widetilde{W}_2$, but by imposing an ``electrical'' condition along the $\vec{e}_3$-direction we exclude the possibility of having certain $\zt$ configurations of $\widetilde{W}_2$ sitting immediately next to each other. These restrictions force $\wel$ to be a proper subsystem of the full $\zz$ extension $\widetilde{W}_2^\zz$. Nevertheless we can show that $\htop(\wel)=\htop(\widetilde{W}_2)= 
h_L(\wel)$.\\

\begin{schritt}[The $\zt$ wire shift]\label{Wstep}
The formal construction of the wire shift $W$ involves an alphabet $\ag_W$ of 7 symbols which we think of as square tiles of unit length as displayed in Figure \ref{Walphfigure}. We will refer to symbol $1$ as the blank symbol and to symbols $2$ up to $7$ containing finite segments of wires (thicklines) as the wire symbols.

\begin{figure}[ht]
\begin{center}
\setlength{\unitlength}{9mm}
\begin{picture}(11.2,2)
\thinlines
\multiput(0,0.5)(1.7,0){7}{\line(0,1){1}}
\multiput(0,0.5)(1.7,0){7}{\line(1,0){1}}
\multiput(1,1.5)(1.7,0){7}{\line(0,-1){1}}
\multiput(1,1.5)(1.7,0){7}{\line(-1,0){1}}
\thicklines
\put(1.7,1){\line(1,0){1}}
\put(3.4,1){\line(1,0){1}}
\put(5.1,1){\line(1,0){1}}
\put(7.3,0.5){\line(0,1){1}}
\put(9,0.5){\line(0,1){1}}
\put(10.7,0.5){\line(0,1){1}}
\put(3.9,1){\line(0,1){0.5}}
\put(5.6,0.5){\line(0,1){0.5}}
\put(9,1){\line(1,0){0.5}}
\put(10.2,1){\line(1,0){0.5}}
\put(3.9,1){\circle*{0.15}}
\put(5.6,1){\circle*{0.15}}
\put(9,1){\circle*{0.15}}
\put(10.7,1){\circle*{0.15}}
\setcounter{help}{1}
\multiput(0.5,0.2)(1.7,0){7}
{\makebox(0,0){\mbox{\small$\arabic{help}$}} \stepcounter{help}}
\end{picture}
\caption{The alphabet $\ag_W$ of the wire shift $W$}
\label{Walphfigure}
\end{center}
\end{figure}
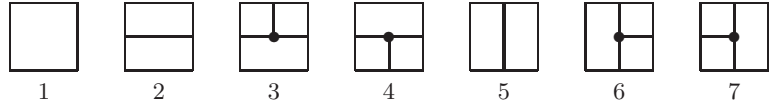

These symbols can be placed next to each other only in a way that conserves wires: Precisely the symbols $2,3,4,7$ with a wire present at their left edge are allowed to sit to the right of a symbol $2,3,4,6$ having a wire at its right edge and exactly the symbols $4,5,6,7$ with a wire present on their lower edge can appear above a symbol $3,5,6,7$ having a wire on its upper edge. Analogously only symbols $1,5,6$ are allowed to the right of $1,5,7$ and one of $1,2,3$ is possible above $1,2,4$. Note that by posing these nearest-neighbor restrictions as adjacency rules we define a non-trivial, even strongly essential (\ie all symbols occur and all allowed transitions are realized in some point) $\zt$ \sft.

\begin{lemma}\label{ExtendibleBlockLemma}
The wire shift $W$ is extendible and block gluing.
\end{lemma}

\begin{proof}
Since $W$ is a nearest-neighbor \sft\ and the horizontal \resp vertical transitions only depend on the wires that are (are not) present at the vertical \resp horizontal borders of each symbol we may consider only the configurations along the boundaries of arbitrary blocks.

Let $B=[\vu,\vv\,]\subsetneq\zt$ be a rectangular block in $\zt$ and $P\in\lang_B(W)$ any allowed configuration on $B$. We define a point $w\in W$ as follows: $w|_B:=P$, and $w|_\vi:=1$ for all $\vi\in\zt\setminus\bigl([\vu-\vec{\idop},\vv+\vec{\idop}\,]\cup (\vu-\vec{\idop}+\zz\vec{e}_1)\cup(\vv+\vec{\idop}+\zz\vec{e}_1)\bigr)$. The remaining coordinates are filled with wire symbols according to Figure \ref{BlockExtendFigure}, where wire segments drawn in grey may or may not be necessary, depending on the symbols along the border of the configuration $P$. It is easily checked that this configuration $w$ truely is a point in $W$. Therefore we have extended an arbitrary valid pattern $P$ on some finite rectangle $B$ to the whole of $\zt$, which proves $W$ to be extendible.

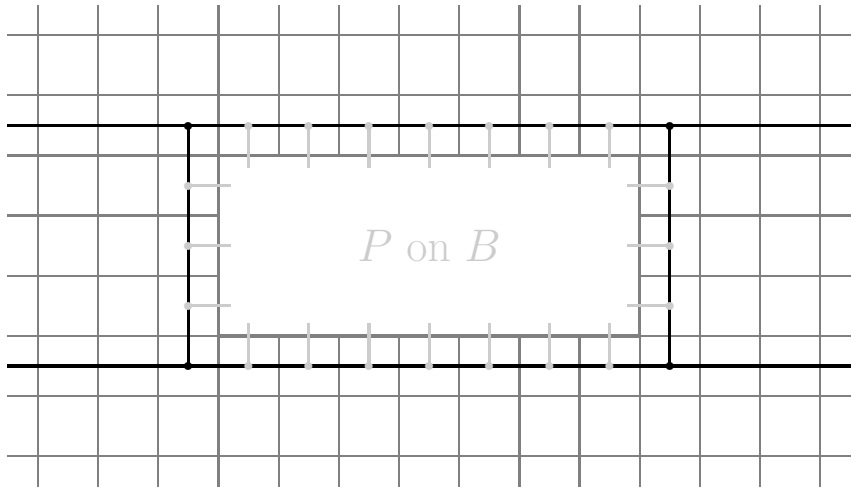
\begin{figure}[ht]
\begin{center}
\setlength{\unitlength}{8mm}
\begin{picture}(14,8)
\thinlines
\color{darkgray}
\multiput(0.5,0)(1,0){4}{\line(0,1){8}}
\multiput(0,0.5)(0,1){3}{\line(1,0){14}}
\multiput(10.5,0)(1,0){4}{\line(0,1){8}}
\multiput(0,5.5)(0,1){3}{\line(1,0){14}}
\multiput(0,3.5)(0,1){2}{\line(1,0){3.5}}
\multiput(10.5,3.5)(0,1){2}{\line(1,0){3.5}}
\multiput(4.5,0)(1,0){6}{\line(0,1){2.5}}
\multiput(4.5,5.5)(1,0){6}{\line(0,1){2.5}}
\thicklines
\put(3.5,2.5){\line(1,0){7}}
\put(3.5,5,5){\line(1,0){7}}
\put(3.5,2.5){\line(0,1){3}}
\put(10.5,2.5){\line(0,1){3}}
\thinlines
\color{black}
\thicklines
\put(0,2){\line(1,0){14}}
\put(0,6){\line(1,0){14}}
\put(3,2){\line(0,1){4}}
\put(11,2){\line(0,1){4}}
\put(3,2){\circle*{0.15}}
\put(3,6){\circle*{0.15}}
\put(11,2){\circle*{0.15}}
\put(11,6){\circle*{0.15}}
\color{gray}
\thicklines
\multiput(4,2)(1,0){7}{\circle*{0.15}}
\multiput(4,6)(1,0){7}{\circle*{0.15}}
\multiput(3,3)(0,1){3}{\circle*{0.15}}
\multiput(11,3)(0,1){3}{\circle*{0.15}}
\multiput(4,2)(1,0){7}{\line(0,1){0.7}}
\multiput(4,6)(1,0){7}{\line(0,-1){0.7}}
\multiput(3,3)(0,1){3}{\line(1,0){0.7}}
\multiput(11,3)(0,1){3}{\line(-1,0){0.7}}
\put(7,4){\makebox(0,0){\mbox{\huge$P$ on $B$}}}
\end{picture}
\caption{Extending finite configurations $P$ on rectangular blocks $B$}
\label{BlockExtendFigure}
\end{center}
\end{figure}

Similarily we may take two patterns $P_1\in\lang_{B_1}(W)$, $P_2\in\lang_{B_2}(W)$ on rectangular blocks $B_1,B_2\subsetneq\zt$ where we suppose $B_1,B_2$ are separated by a distance larger than $2$. There are two cases: Either $B_1$, $B_2$ are separated by a distance $>2$ along direction $\vec{e}_2$ (suppose this is the vertical direction in Figure \ref{BlockExtendFigure}) or along direction $\vec{e}_1$. In the first case we surround each pattern $P_i$ ($i=1,2$) with a wire as in Figure \ref{BlockExtendFigure}.  Note that the $\vec{e}_2$-separation of at least two coordinates is large enough to do this without causing any conflict in placing symbols. Hence -- filling all remaining coordinates with blanks -- there is a valid point $w\in W$ realizing both patterns, \ie $w|_{B_1}=P_1$ and $w|_{B_2}=P_2$. For the second case, $B_1,B_2$ being separated along direction $\vec{e}_1$ by a distance $>2$, we just have to rotate the picture in Figure \ref{BlockExtendFigure} by $90$ degrees, which can be done as the alphabet and the transition rules are invariant under symbol rotation, and proceed as before. This proves $W$ being block gluing (at gap $g=2$) as claimed.
\end{proof}

\begin{rem}
Since the wire shift $W$ is block gluing and has an alphabet with more than one symbol its topological entropy is strictly positive. Calculations give the estimate $\log 1.75<\htop(W)<\log 1.964$. Moreover every non-trivial continuous factor of $W$ is block gluing again and thus has to have strictly positive entropy, \ie $W$ has topologically completely positive entropy.
\end{rem}
\end{schritt}

\begin{schritt}[Splitting the blank symbol]\label{Wsplittstep}
Next we modify $W$ by splitting the blank symbol into $k\geq 2$ distinct, but completely interchangeable copies, thus the new alphabet $\ag_k:=\{1_i\mid 1\leq i\leq k\}\,\dcup\,\{2,3,4,5,6,7\}$ has $(k+6)$ symbols. The adjacency rules concerning wire-conservation stay unchanged and it can be checked that the above argument showing extendability and the block gluing property for $W$ (Lemma \ref{ExtendibleBlockLemma}) does not at all depend on the number of distinct types of blank symbols and thus just carries over. Hence we immediately get

\begin{corol}
For $k\geq 2$ the new $\zt$ \sft s $\widetilde{W}_k$ again are extendible and block gluing (at gap $g=2$).
\end{corol}

In addition we are able to calculate the topological entropy and show that $\widetilde{W}_k$ is no longer entropy-minimal.

\begin{lemma}\label{EntropyLemma}
For every $k\geq 2$ the topological entropy of $\widetilde{W}_k$ equals $\log k$. In particular each $\widetilde{W}_k$ contains a proper $\zt$ \sft\ subsystem of full entropy.
\end{lemma}

\begin{proof}
As there are no restrictions on the types of blank symbols that can be placed next to each other horizontally or vertically, $\widetilde{W}_k$ contains the full shift on its $k$ blanks as a subsystem and thus $\htop(\widetilde{W}_k)\geq \log k$.

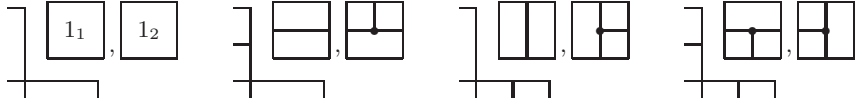
\begin{figure}[ht]
\begin{center}
\setlength{\unitlength}{7.5mm}
\begin{picture}(15.4,2)
\thinlines
\multiput(0.2,0.5)(4,0){4}{\line(1,0){1.6}}
\multiput(0.5,0.2)(4,0){4}{\line(0,1){1.6}}
\multiput(1.8,0.2)(4,0){4}{\line(0,1){0.3}}
\multiput(0.2,1.8)(4,0){4}{\line(1,0){0.3}}
\multiput(0.9,0.9)(4,0){4}{\line(0,1){1}}
\multiput(0.9,0.9)(4,0){4}{\line(1,0){1}}
\multiput(1.9,0.9)(4,0){4}{\line(0,1){1}}
\multiput(0.9,1.9)(4,0){4}{\line(1,0){1}}
\multiput(2.2,0.9)(4,0){4}{\line(0,1){1}}
\multiput(2.2,0.9)(4,0){4}{\line(1,0){1}}
\multiput(3.2,0.9)(4,0){4}{\line(0,1){1}}
\multiput(2.2,1.9)(4,0){4}{\line(1,0){1}}
\thicklines
\put(4.2,1.15){\line(1,0){0.3}}
\put(9.15,0.2){\line(0,1){0.3}}
\put(12.2,1.15){\line(1,0){0.3}}
\put(13.15,0.2){\line(0,1){0.3}}
\put(4.9,1.4){\line(1,0){1}}
\put(6.2,1.4){\line(1,0){1}}
\put(12.9,1.4){\line(1,0){1}}
\put(9.4,0.9){\line(0,1){1}}
\put(10.7,0.9){\line(0,1){1}}
\put(14.7,0.9){\line(0,1){1}}
\put(6.7,1.4){\line(0,1){0.5}}
\put(10.7,1.4){\line(1,0){0.5}}
\put(13.4,1.4){\line(0,-1){0.5}}
\put(14.7,1.4){\line(-1,0){0.5}}
\put(6.7,1.4){\circle*{0.15}}
\put(10.7,1.4){\circle*{0.15}}
\put(13.4,1.4){\circle*{0.15}}
\put(14.7,1.4){\circle*{0.15}}
\setcounter{help}{1}
\multiput(1.4,1.4)(1.3,0){2}
{\makebox(0,0){\mbox{\small$1_\arabic{help}$}} \stepcounter{help}}
\multiput(2.05,1)(4,0){4}
{\makebox(0,0){\mbox{\small$,$}}}
\end{picture}
\caption{$\widetilde{W}_2$ has corner condition 2}
\label{cornercondFigure}
\end{center}
\end{figure}

To show the reversed inequality we look at the case $k=2$ and check that $\widetilde{W}_2$ actually has corner condition $2$, \ie for any given configuration on the coordinates $\vi+\{(-1,-1),(-1,0),(0,-1)\}$ there are exactly two symbols possible at coordinate $\vi\in\zt$ (see Figure \ref{cornercondFigure}). For $k>2$, given a corner configuration of the first type (left-most) in Figure \ref{cornercondFigure} filling in one of the blanks yields precisely $k$ possibilities, whereas for corner configurations of the other three types the number of possible symbols stays $2$. Hence none of the cases give more than $k$ choices for the next symbol. Following from this any configuration on the lower left border $R:=\{\vi\in C_n\mid \vi_1=-n\,\vee\,\vi_2=-n\}$ of a big square $C_n:=[-n\vec{\idop},n\vec{\idop}]$ for $n\in\nz$ allows for at most $k^{4n^2}$ ways to fill in the remaining $4n^2$ coordinates of $C_n$. Summing over all configurations on $R$ we get a coarse estimate
\[
\big\vert\lang_{C_n}(\widetilde{W}_k)\big\vert\leq \sum_{\lang_{R}(\widetilde{W}_k)} k^{4n^2}\leq \big\vert\ag_k\big\vert^{\abs{R}}\cdot k^{4n^2}=(k+6)^{4n+1}\cdot k^{4n^2}\ .
\]
Putting this bound into the definition of topological entropy gives
\[
\htop(\widetilde{W}_k)= \lim_{n\rightarrow\infty} \frac{\log\big\vert\lang_{C_n}(\widetilde{W}_k)\big\vert}{\big\vert C_n\big\vert}\leq \lim_{n\rightarrow\infty}\frac{\log\big((k+6)^{4n+1}\cdot k^{4n^2}\big)}{4n^2+4n+1}=\log k\ .
\]
So $\htop(\widetilde{W}_k)=\log k$ as claimed and the full shift on the $k$ blank symbols constitutes a proper \sft\ subsystem of full entropy.
\end{proof}

\begin{rem}
We do not know whether the wire shift $W$ from Step \ref{Wstep} itself is entropy minimal or not. However the above shows that there is a conceptual change in entropy between the case of only one blank in $W$ with $\htop(W)>\log 1$ and the case of $k\geq 2$ blanks in $\widetilde{W}_k$ with $\htop(\widetilde{W}_k)=\log k$. So having (at least) the two blanks is enough to compensate the entropy of the underlying shift $W$, thus generating a proper subsystem of equal entropy. It seems that this phenomenon does occur in many $\zd$ shifts (see \eg the meandering streams shift in \cite{bps}) where we can eventually get non entropy minimal systems by independently splitting an apropriate symbol (in the presence of a fixed point) or a couple of symbols (in the presence of some periodic point).
\end{rem}

To see that none of the \sft s $\widetilde{W}_k$ is actually topologically conjugate to a full shift it suffices to count the fixed points of $\widetilde{W}_k$ -- there are $k+2$.

\begin{lemma}\label{WknondegLemma}
None of the modified wire shifts $\widetilde{W}_k$ ($k\geq 2$) is degenerate with respect to any sublattice. The same is true for $W$.
\end{lemma}

\begin{proof}
Let $L=\seq{\vu}_\zz\subsetneq\zt$ be any 1-dimensional sublattice; then at least one of the two standard base vectors $\vec{e}_1,\vec{e}_2\in\zt$ does not belong to $L$. Fixing a complementary sublattice $L'=\seq{\vv}_\zz\subsetneq\zt$ we may assume that $\vec{e}_1=m\vu+n\vv$ (\resp $\vec{e}_2=m\vu+n\vv$) with $m,n\in\zz$ and $n\neq 0$. As $\bigl(\widetilde{W}_k\bigr)_L$ contains points seeing an arbitrary symbol of $\ag_k$ at a particular coordinate in $L$ we may pick $w^{(1)},w^{(2)}\in\bigl(\widetilde{W}_k\bigr)_L$ such that $w^{(1)}_{\vec{0}}=7$ and $w^{(2)}_{m\vu}=2$. Now if $\widetilde{W}_k$ were degenerate with respect to $L$, \ie $\widetilde{W}_k=\bigl\{\bigl(w^{(\vi)}\in\bigl(\widetilde{W}_k\bigr)_L\bigr)_{\vi\in L'}\bigr\}$, a family $\bigl(w^{(\vi)}\bigr)_{\vi\in L'}$ with $w^{(\vec{0})}:=w^{(1)}$ and $w^{(n\vv)}:=w^{(2)}$ would give a valid point $w:=\bigl(w^{(\vi)}\bigr)_{\vi\in L'}\in \widetilde{W}_k$. However this point $w$ would see a symbol $7$ at the origin and a symbol $2$ at coordinate $\vec{e}_1$ (\resp $\vec{e}_2$), which is not possible due to the transition rules forcing wire conservation.
\end{proof}

Looking at the projectional entropies of $\widetilde{W}_k$ we get that the infimum of $h_L(\widetilde{W}_k)$ taken over all sublattices $L\subsetneq\zt$ is attained exactly in the two principal directions. Moreover this infimum is strictly bounded away from the topological entropy $\htop(\widetilde{W}_k)$.

\begin{lemma}
The projectional entropy of $\widetilde{W}_k$ ($k\geq 2$) with respect to the horizontal respectively vertical sublattice $L_1=\seq{\vec{e}_1}_\zz\subsetneq\zt$ \resp $L_2=\seq{\vec{e}_2}_\zz\subsetneq\zt$ is $h_{L_1}(\widetilde{W}_k)=h_{L_2}(\widetilde{W}_k)=\log\bigl(2+\frac{k}{2}+\frac{\sqrt{k^2-4k+8}}{2}\bigr)$. For every other sublattice $L\subsetneq\zt$ the projectional entropy is $h_L(\widetilde{W}_k)=\log(k+6)$.
\end{lemma}

\begin{proof}
For $L_1$ the projection $\bigl(\widetilde{W}_k\bigr)_{L_1}$ is the 1-dimensional subshift given by the nearest neighbor transition conditions that ensure wire conservation in the horizontal direction, \ie a sequence $(w_i)_{i\in\zz}\in{\ag_k}^\zz$ of symbols is valid, if and only if the wire is conserved in each subword $w_i\,w_{i+1}$ ($i\in\zz$). Thus the projectional entropy is given as the logarithm of the Perron-eigenvalue of the corresponding transition matrix. A simple calculation gives $\lambda_{\text{P}}=2+\frac{k}{2}+\frac{\sqrt{k^2-4k+8}}{2}$.

As the definition of the wire shifts $\widetilde{W}_k$ is completely symmetric with respect to the $\vec{e}_1$- and $\vec{e}_2$-direction the same holds for the sublattice $L_2$.

For other sublattices $L\subsetneq\zt$ different coordinates in $L$ are never horizontally or vertically adjacent. Thus we can freely place symbols on all coordinates in $L$. It is then easily checked that any configuration on $L$ can be extended obeying the rules on wire conservation to get a configuration on all of $\zt$ (recall Figure \ref{cornercondFigure} to see that every corner can be filled). This gives a valid point of $\widetilde{W}_k$. Hence $\bigl(\widetilde{W}_k\bigr)_L={\ag_k}^\zz$ and $h_L(\widetilde{W}_k)=\log\bigl|\ag_k\bigr|$.
\end{proof}
\end{schritt}

In the following we are going to build a $\zz^3$ \sft\ and for this we will take the $\vec{e}_3$-direction as vertical, whereas the $\vec{e}_1$- and $\vec{e}_2$-axes are situated in a horizontal plane. Therefore from now on ``above'' resp.\ ``below'' refers to an increase resp.\ decrease of the $\vec{e}_3$-component of $\zz^3$ coordinates.

\begin{schritt}[The $\zz^3$ electrical wire shift]\label{Welstep}
Let $L:=\seq{\vec{e}_1,\vec{e}_2}_\zz\subsetneq\zz^3$. To get our $\zz^3$ \sft\ $\wel$ we stockpile infinitely many $\zt$ configurations of $\widetilde{W}_2$ as horizontal layers in a point in $\wel$, \ie $\wel\subseteq\bigl\{w\in{\ag_2}^{\zz^3}\bigm\vert \forall\,n\in\zz:\ w|_{L+n\vec{e}_3}\in\widetilde{W}_2\bigr\}=\widetilde{W}_2^\zz$. In order to get $\wel$ being a proper subset of $\widetilde{W}_2^\zz$ we pose the following ``electrical'' restriction on the allowed vertical transitions: Whenever a wire symbol $2$ resp.\ $5$ appears at a coordinate $\vi\in\zz^3$ then at coordinates $\vi\pm\vec{e}_3$ we are only allowed to see an arbitrary blank or a symbol $5$ resp.\ $2$. If we have a symbol $a\in\{3,4,6,7\}$ at some coordinate then there cannot be any wire symbol directly above or below; instead we have to see one of the blanks there. We can think of this condition as forbidding parallel wires next to each other in the $\vec{e}_3$-direction, as they may cause ``interference'' of signals. However note that wires may and will cross on vertically adjacent horizontal layers (placing a symbol $2$ above a symbol $5$ or vice versa). As this extra condition is still given by nearest-neighbor restrictions, $\wel\subsetneq\widetilde{W}_2^\zz$ is a \sft.

The projection map $\pi_L:\;\wel\rightarrow\widetilde{W}_2,\ w\mapsto w|_L$ restricting points in $\wel$ to the horizontal sublattice $L$ is surjective (taking any $\widetilde{W}_2$-point as a configuration on $L$ and filling $\zz^3\setminus L$ with blanks gives rise to a valid point in $\wel$), so $\wel_{L}=\widetilde{W}_2$ and then $\wel\subsetneq (\wel_L)^\zz$.

\begin{lemma}
$\wel$ is extendible and block gluing.
\end{lemma}

\begin{proof}
The argument is a slight elaboration of the proof we presented for Lemma \ref{ExtendibleBlockLemma}. Taking any pattern $P\in\lang_B(\wel)$ on some finite 3-dimensional cuboid $B:=[\vu,\vv\,]\subsetneq\zz^3$, to build a valid point $w\in\wel$ with $w|_B=P$, we surround every other (finite) horizontal layer $B\cap(L+n\vec{e}_3)$ (with $n\in 2\zz\cap[\vu_3,\vv_3]$) of $P$ with a wire exactly as in Figure \ref{BlockExtendFigure}. For each remaining horizontal layer $H':=B\cap(L+n\vec{e}_3)$ with $n\in (2\zz+1)\cap[\vu_3,\vv_3]$ fixed, we first enlarge the rectangular pattern $P|_{H'}$ to a pattern on the larger rectangle $B':=H'+[-\vec{\idop},\vec{\idop}\,]\subsetneq\zt$ (adding a border of width 1 to $H'$) by extending all wires hitting the boundary of $P|_{H'}$ by another straight segment, \ie a symbol $2$ or $5$. If at a certain coordinate no wire segment hits the boundary of $H'$, we fill the adjacent coordinate in $B'\setminus H'$ with a blank. The four corners of $B'$ are also filled with (arbitrary) blanks. After this step we still have a valid pattern $P'\in\lang_{B'}(\widetilde{W}_2)$ with $P'|_{H'}=P|_{H'}$. As the pattern $P$ is valid for $\wel$ we know that above and below any of the newly placed wire symbols in $B'\setminus H'$ there can be no symbols $3,4,6,7$, as otherwise the ``electrical'' restrictions on vertical adjacencies would be violated already in $P$. By construction of the surrounding wires on the horizontal layers with even $\vec{e}_3$-coordinates in the previous step, above and below a symbol $2$ in the border of $B'$ we always see a symbol $5$ and above and below a symbol $5$ there always is a symbol $2$. As those transitions are allowed vertically (crossing wires do not cause electrical interference), we still have a valid configuration for $\wel$. Now we can surround the enlarged rectangular pattern $P'$ by a wire and fill all remaining coordinates in the horizontal layer $L+n\vec{e}_3$ with blanks. Again there is no conflict at the four coordinates where the wire surrounding $P'$ crosses the wires already put to surround our pattern $P$ in the horizontal layers $L+(n\pm 1)\vec{e}_3$. As this can be done independently for all choices of $n\in (2\zz+1)\cap[\vu_3,\vv_3]$ we have constructed a valid configuration for $\wel$ on all of $B+L\subsetneq\zz^3$. Finally we may fill the coordinates in $\zz^3\setminus(B+L)$ with blanks to obtain a valid point in $\wel$.

A similar reasoning shows that $\wel$ is block gluing. Suppose we are given two finite cuboids $B_1,B_2\subsetneq\zz^3$ which are separated along the $\vec{e}_3$-direction by at least distance $2$. Extending arbitrary valid patterns on $B_1,B_2$ can be done as in the last paragraph and filling in the horizontal layer(s) between them as well as all the remaining coordinates with blanks will produce a valid point in $\wel$ which realizes both patterns. For $B_1$, $B_2$ having a distance larger than $4$ along the $\vec{e}_2$-direction we can do the same. Note that wiring all horizontal layers of our patterns on $B_1,B_2$ does not need a space of more than $2$ coordinates in direction $\vec{e}_2$, thus a separation of $4$ symbols is sufficient for not causing any conflict. The case of $\vec{e}_1$-separation by a distance larger than $4$ is then immediate -- just recall the possibility of rotating Figure \ref{BlockExtendFigure} by $90$ degrees. Therefore our $\zz^3$ \sft\ $\wel$ is block gluing (at gap $g=4$).
\end{proof}

Since $\wel$ contains 14 coordinatewise periodic points of period $(1,1,2)$, it is not conjugate to a full shift (which would contain $\abs{\ag}^2$ of those points).

As before we can determine the precise value of the topological entropy: $\wel$ being a subset of $\widetilde{W}_2^\zz$ forces $\htop(\wel)\leq \htop(\widetilde{W}_2^\zz)=\log 2$ whereas $\wel$ still containing the full shift on 2 (blank) symbols implies $\htop(\wel)\geq\log 2$, thus $\htop(\wel)=\log 2=h_L(\wel)$ and we have verified all the properties we claimed in Proposition \ref{WelCounterexampleprop}.\\

Moreover we have the following:

\begin{lemma}
The electrical wire shift $\wel$ is not degenerate with respect to any sublattice.
\end{lemma}

\begin{proof}
The argument is a simple generalization of the proof used for Lemma \ref{WknondegLemma}.
\end{proof}

\begin{questions}
What is the projectional entropy of $\wel$ with respect to any 1-dimensional sublattice? With respect to any 2-dimensional one? (The infimum equals exactly the topological entropy and is attained at least for the horizontal sublattice $L=\seq{\vec{e}_1,\vec{e}_2}\subsetneq\zz^3$.)

Is any conjugate presentation of $\widetilde{W}_k$ or $\wel$ degenerate with respect to any sublattice $L$ (I would guess not)?

Can we (inductively) calculate the number of coordinatewise periodic points of $\widetilde{W}_k$ or $\wel$ in order to show that there is no degeneracy? (Using the transition matrices, which have a ``nice'' recursive structure.)\\
What about the image of the full shift on $k$ blanks under any conjugacion (it is again a proper subsystem)? Can we use this to show there is no degeneracy?
\end{questions}
\end{schritt}

\begin{ack}
Part of this paper was written during a short stay at Erwin Schr\"odinger Institute, Vienna and the author would like to thank Prof.\ Klaus Schmidt for his invitation and the kind hospitality provided by his institute.
\end{ack}


\begin{thebibliography}{99}
\bibitem{bps}
M. Boyle, R. Pavlov and M. Schraudner, {\em Multidimensional sofic shifts without separation, and their factors}, preprint.
\bibitem{dgs}
M. Denker, C. Grillenberger and K. Sigmund, {\em Ergodic theory on compact spaces}, Lecture Notes in Mathematics 527, Springer, Berlin (1976).
\bibitem{jkm}
A. Johnson, S. Kass and K. Madden, {\em Projectional entropy in higher dimensional shifts of finite type}, Complex Systems {\bf 17} (2007), 243--257.
\bibitem{lm}
D. Lind and B. Marcus, {\em Introduction to symbolic dynamics and coding.} Cambridge University Press, Cambridge (1995).
\bibitem{ppc}
R.\ Pavlov, personal communication (2008).
\bibitem{qs}
A. Quas and A. Sahin, {\em Entropy gaps and locally maximal entropy in $\zd$ subshifts}, Ergodic Theory Dynam. Systems {\bf 23} (2003), 1227-1245.
\bibitem{rs}
E.A. Robinson and A. Sahin, {\em Mixing properties of nearly maximal entropy measures for $\zd$ shifts of finite type}, Colloq. Math. {\bf 84/85} (2000), part 1, 43--50.
\end{thebibliography}
\end{document}